\input amstex 
\hsize=30truecc
\documentstyle{amsppt} 
\pagewidth{30pc}
\pageheight{52pc}
\loadbold
\TagsOnRight
\nologo
\nopagenumbers

\define\na{\Bbb N} 
\define\proc{$(\Cal P^\z,\mu,\sigma)$}

\define\procbr{$((\Cal P_B^r)^\z,\mu_B,\sigma_B)$}
 
\define\sq{sequence} 

\define\z{\Bbb Z}
\define\r{\Bbb R}

\topmatter  
 
\title The law of series\endtitle
\author 
T. Downarowicz and Y. Lacroix
\endauthor 
\address 
Institute of Mathematics and Computer Science, Wroclaw University of
Technology, Wy\-brze$\dot{\text z}$e Wys\-pia{\'n}\-skie\-go 27, 
50@-370 Wroc{\l}aw, 
Poland
\endaddress
\date January 8, 2006\enddate
\email 
downar\@pwr.wroc.pl
\endemail
\thanks
This paper was written during the first author's visit at CPT/ISITV, supported by
CNRS. The research of the first author is supported by DBN grant 1 P03A
021 29  
\endthanks 
\address 
Institut des Sciences de l'Ing\'enieur de Toulon et du Var,
Avenue G. Pompidou, B.P. 56, 83162 La Valette du Var Cedex,
France
\endaddress
\email 
yves.lacroix\@univ-tln.fr
\endemail

\abstract 
We prove a general ergodic-theoretic result concerning the return time statistic, which, properly
understood,  sheds some new light on the common sense phenomenon known as
{\it the law of series}. Let \proc\ be an ergodic 
process on finitely many states, with positive entropy. We show that the distribution function of the
normalized  waiting time for the first visit to a small cylinder set $B$
is, for majority of such cylinders and up to epsilon,  
dominated by the exponential distribution function $1-e^{-t}$. This fact has the following
interpretation:  The occurrences of such a ``rare event'' $B$ can deviate
from purely random in only one direction -- so that for 
any length of an ``observation period'' of time, the first occurrence of $B$ ``attracts'' its further
repetitions  in this period.
\endabstract 

\endtopmatter 
\document

\heading Note \endheading
This paper resulted from studying asymptotic laws for return/hitting time statistics in stationary
processes, a field in ergodic theory rapidly developing in the recent years (see e.g. [A-G], [C], [C-K], 
[D-M], [H-L-V], [L] and the reference therein). Our result significantly contributes to this area due to 
both its generality
and strength of the assertion. After having completely written the proof, during a free-minded discussion, 
the authors have discovered an astonishing interpretation of the result, clear even in terms of the common 
sense understanding of random processes. The organization of the paper is aimed to emphasise this discovery. 
The consequences for the field of asymptotic laws are moved toward the end of the paper.

\heading Introduction\endheading 
The phenomenon known as {\it the law of series} appears in many aspects of every-day life, science and technology. 
In the common sense understanding it signifies a sudden increase of frequency of a rare event, seemingly violating 
the rules of probability. Let us quote from Jean Moisset ([Mo]):
\block{\eightpoint This law can be defined as: the repetition of identical or analogous events, things or symbols in space or time; 
for example: the announcement of several similar accidents on the same day, a series of strange events experienced by someone on the same day which are either
happy ones (a period of good luck) or unfavorable (disastrous) ones, or the repetition of unexpected similar events. For example, you are invited to dinner and you
are served a roast beef and you note that you were served the same menu the day before at your uncle's home and the day before that at your cousin's home. }
\endblock
Another proverb describing more or less the same (with regard to unwanted events) is {\it misfortune seldom 
comes alone}. Both expressions exist in many languages, proving that the phenomenon has been commonly noticed 
throughout the world.
In this setting it has been accounted to the category of unexplained mystery, paraphysics, parapsychology, 
together with ``malignancy of fate'', ``Murphy laws'', etc. Many pseudoscientific experiments have been 
conducted to prove obvious violation of statistical laws, where such laws were usually identified with 
the statistics of a series of independent trials ([St], [Km]). Equally many texts have been devoted to 
explain the anomaly within the framework of an independent process (see e.g., [Mi]), or merely as the weakness 
of our memory, keen to notice unusual events as more frequent just because they are more distinctive.

The phenomenon is also known in more serious science. In modeling and statistics it is sometimes
called ``clustering of data''. It is experimentally observed in many real processes, such as traffic jams, 
telecommunication network overloads, power consumption peaks, demographic peaks, stock market fluctuations, 
etc., as periods of increased frequency of occurrences of certain rare events. These anomalies are usually 
explained in terms of physical dependence (periods of propitious conditions) and complicated algorithms 
are implemented in modeling these processes to simulate them. 

But, to our knowledge, there was no logical construction proving, in full generality, that there 
exists a  ``natural'' tendency of rare events to appear in series, and the result we will present 
in this paper, or any of its possible variants, remained until now unnoticed by the specialists. 
We prove an ergodic-theoretic theorem on stationary stochastic processes, in which a wide range 
of ``rare events'' is shown to behave either in a way which we call ``unbiased'' (i.e., as in an
independent process), or else exactly as it is specified in {\it the law of series}, i.e., so that 
the first occurrence increases the chances of untimely repetitions. Roughly speaking, we prove that 
{\bf rare events appear in series whenever the unbiased behavior is perturbed: there is no other choice}. 
Besides ergodicity (which is automatically satisfied if we observe a chosen at random single realization 
of any process) we make only one essential, but obviously necessary, assumption on the process: it must 
maintain a ``touch of randomness'', i.e., the future must not be completely determined by the past, 
which is equivalent to assuming positive entropy. Without this assumption a rotation of a compact group 
is an immediate example where events never appear in series. Of course, not every interesting rare event 
in reality can be modeled by the type of set we describe (cylinder over a long block), and we do not claim 
that our theorem fully explains the common sense phenomenon, but it certainly sheds on it some new light. 

In terms of ergodic theory, we define two elementary antagonistic properties of the return times
called ``attracting'' and ``repelling'', and we prove that they behave quite differently in processes of zero 
and of positive entropy: attracting can persist for arbitrarily long blocks in both cases, while repelling 
must decay (as the length of blocks grows to infinity) in positive entropy processes. 
Many properties are known to differentiate between positive and zero entropy, but most of them involve a 
passage via measure-theoretic isomorphism, i.e., change of a generator, or require some additional structure. 
Our ``decay of repelling'' holds in general and for any finite generator, or even partition, as long as it 
generates positive entropy. 

It is impossible not to mention here the theorem of Ornstein and Weiss [O-W2] which relates the return 
times of long blocks to entropy. However, this theorem says nothing about attracting or repelling,
because the limit appearing in the statement is insensitive to the proportions between gap sizes. 
Nevertheless, this remarkable result is very useful and it will help also in our proof. 
Our theorem's proof is entirely contained within the classics of ergodic theory; it relies 
on basic facts on entropy for partitions and sigma-fields, some elements of the Ornstein theory 
($\epsilon$-independence), the Shannon-McMillan-Breiman Theorem, the Ornstein-Weiss Theorem on return times, 
the Ergodic Theorem, basics of probability and calculus. We do not invoke any specialized machinery of 
stochastic processes or statistics.
\medskip
The authors would like to thank Dan Rudolph for a hint leading to the construction of Example 2 and, in effect,
to the discovery of the attracting/repelling asymmetry. We also thank Jean-Paul Thouvenot for his interest in
the subject, substantial help, and the challenge to find a purely combinatorial proof (which we save for 
the future).

\heading Rigorous definition and statement \endheading
We establish the notation necessary to formulate the main result. Let \proc\ be an ergodic process on 
finitely many symbols, i.e., $\#\Cal P<\infty$, $\sigma$ is the standard left shift map and $\mu$ is an 
ergodic shift-invariant probability measure on $\Cal P^\z$. Most of the time, we will identify finite blocks 
with their 
cylinder sets, i.e., we agree that $\Cal P^n=\bigvee_{i=0}^{n-1}\sigma^{-i}(\Cal P)$. Depending on the context, 
a block $B\in\Cal P^n$ is attached to some coordinates or it represents a ``word'' which may appear in different 
places along the $\Cal P$-names. We will also use the probabilistic language of random variables. Then 
$\mu\{R\in A\}$ ($A\subset\r$) will abbreviate $\mu(\{x\in\Cal P^\z:R(x)\in A\})$. Recall, that if the random 
variable $R$ is nonnegative and $F(t)=\mu\{R\le t\}$ is its distribution function, then the expected value 
of $R$ equals $\int_0^\infty 1-F(t)\,dt$.

For a set $B$ of positive measure let $R_B$ and $\overline R_B$ denote the random variables defined on $B$ 
(with the conditional measure $\mu_B =\frac\mu{\mu(B)}$) as the absolute and normalized first return time to $B$, 
respectively, i.e.,
$$
R_B(y) = \min\{i>0, \sigma^i(y)\in B\}, \ \ \overline R_B(y) = \mu(B)R_B(y).
$$
Notice that, by the Kac Theorem ([Kc]), the expected value of $R_B$ equals $\frac1{\mu(B)}$, hence that of 
$\overline R_B$ is 1 (that is why we call it ``normalized''). We also define
$$
G_B(t) = \int_0^t 1-F_B(s)\, ds.
$$
(The interpretation of this function is discussed in the following section.)
Clearly, $G_B(t)\le \min\{t,1\}$ and the equality holds when $F_B(t) = 1_{[1,\infty)}$, that is, when $B$ 
occurs precisely with equal gaps, i.e., periodically; the gap size then equals $\frac1{\mu(B)}$.

The key notions of this work are defined below:
\definition{Definition 1}
We say that the visits to $B$ {\it repel} (resp\. {\it attract}) each other with intensity $\epsilon$ from a 
distance $t>0$, if 
$$
G_B(t)\ge 1-e^{-t}+\epsilon\ \ \  (\text{resp\. if }G_B(t)\le 1-e^{-t}-\epsilon).
$$
We abbreviate that $B$ {\it repels} ({\it attracts}) with intensity $\epsilon$ if its visits repel (attract)
each other with intensity $\epsilon$ from some distance $t$.
\enddefinition
Obviously, occurrences of an event may simultaneously repel from one distance and attract from
another. Notice, that the maximal intensity of repelling is $e^{-1}$ achieved at $t=1$ when $B$ 
appears periodically. The intensity of attracting can be arbitrarily close to 1 (when $B$ appears 
in enormous clusters separated by huge pauses; see the next section). The main result follows:

\proclaim{Theorem 1} If \proc\ is ergodic and has positive entropy, then for every \hbox{$\epsilon>0$} 
the measure of the union of all $n$-blocks $B\in\Cal P^n$ which repel with intensity $\epsilon$, converges 
to zero as $n$ grows to infinity.
\endproclaim

We also provide an example (Example 2) in which, for a substantial collection of lengths, the majority of 
cylinders display strong attracting. Moreover, the process of Example 2 is isomorphic to a Bernoulli 
process, which implies that a partition with such strong attracting properties can be found in any 
measure-preserving transformation of positive entropy (see the Remark 3). Let us mention that it is 
easy to find zero entropy examples with either persistent repelling (discrete spectrum) or attracting 
(see the Example 3), or even both at a time (see the Remark 4).

\heading Interpretation and its limits\endheading
Let us elaborate a bit on the meaning of attracting and repelling for an event $B$. Let $V_B$ be the random
variable defined on $X$ as the {\it hitting time statistic}, i.e., the waiting time for the first visit in 
$B$ (the defining formula is the same as for $R_B$, but this time it is regarded on $X$ with the measure $\mu$). 
Further, let $\overline V_B = \mu(B)V_B$, called, by analogy, {\it the normalized hitting time} (although 
the expected value of this variable need not be equal to 1). By ergodicity, $V_B$ and $\overline V_B$ are well 
defined. By an elementary consideration of the skyscraper above $B$, one easily verifies, that the 
distribution function $\tilde F_B$ of $\overline V_B$ satisfies the inequalities: 
$$
G_B(t) -\mu(B)\le \tilde F_B(t) \le G_B(t)
$$
(see [H-L-V] for more details). 
Because we deal with long blocks (so that, by the Shannon-McMillan-Breiman Theorem, $\mu(B)$ is, with high 
probability, very small), for sake of the interpretation, we will simply assume that $\tilde F_B = G_B$. 
Thus, attracting and repelling can be considered properties of the hitting rather than return time statistic. 
In fact, if we replace $G_B$ by $\tilde F_B$ in the definition of attracting/repelling, the formulation of 
Theorem 1 remains exactly the same, because it admits tolerance up to a fixed $\epsilon$.

It is easy to see that if \proc\ is an independent Bernoulli process, then, for any long block $B$, 
$F_B(t) \approx 1-e^{-t}$ (and also $\tilde F_B(t) \approx 1-e^{-t}$) with high uniform accuracy. We will 
call such behavior ``unbiased'' (neither attracting nor repelling), and attracting and repelling can be viewed 
as deviations from the unbiased pattern.

Fix some $t>0$. Consider the random variable $I$ counting the number of occurrences of $B$ in the time period \hbox{$[0, \frac t{\mu(B)}]$}. The expected value of $I$ equals $\mu(B)\lfloor\frac t{\mu(B)}\rfloor\approx t$ 
(up to the ignorable error $\mu(B)$). On the other hand, $\mu\{I>0\} = \mu\{V_B\le \frac t{\mu(B)}\} = \tilde F_B(t)$. 
Attracting from the distance $t$ occurs when the last value is smaller than it would be (say, for the same cylinder 
$B$) in an independent process. Because the expected value of $I$ in the independent process is maintained 
(and equals approximately $t$), the conditional expected value of $I$ on the set $\{I>0\}$ must be larger in \proc\ 
than in the independent process. This fact can be interpreted as follows:  If we observe the process for time 
$\frac t{\mu(B)}$ (which is our ``memory length'' or ``lifetime of the observer'') and we happen to see the event 
$B$ during this time at least once, then the expected number of times we will observe the event $B$ is 
larger than the analogous value in the independent process. The first occurrence of $B$ ``attracts'' further 
repetitions ({\it misfortune seldom comes alone}). 
$$
\align 
\text{repelling \ }& ..{...B}......\underline{B......B...B.....B}...{...B}.....
\underline{B....B}.....\underline{B....B..B}....B...{...B}..\\
\text{unbiased \ }& 
..{...B}........\underline{B....B..B....B}.....{...B}......\underline{B..B}
.......\underline{B..B.B}...{...B}...{...B}..\\
\text{attracting \ }&..{...B}..........\underline{B..B.B..B}........
{...B}.......\underline{BB}.........\underline{B.BB}....{...B}...{...B}..\\
\vspace {5pt}
\text{strong attr. \ }
&.........BBB.BB.....................B.................BBB.BB.BB.............\\
\endalign
$$ 
{\eightpoint\it Figure 1: Comparison between unbiased, repelling and attracting distributions of copies of a block.
Attracting with intensity close to 1 occurs, when $G_B$ is very ``flat'' (close to zero on a long initial interval). 
Then $F_B$ is immediately very close to 1 indicating that on most of $B$ the first return time is much smaller 
than $\frac 1{\mu(B)}$. Of course, this must be compensated on a small part of $B$ by extremely large values
of the return time. This means that the visits to $B$ occur in enormous clusters of very high frequency, compensated 
by huge pauses with no (or very few) visits. Such pattern will be called ``strong attracting'' and it will take 
place in some of our examples.} 

\medskip\noindent
Repelling from the distance $t$ means exactly the opposite: The first occurrence lowers the expected number of 
repetitions within the observation period, i.e., repels them. If we have a mixed behavior, our impression about 
whether the event attracts or repels its repetitions depends on the length of our ``memory''. Attracting not assisted 
by repelling (or assisted by repelling of an ignorably small intensity) means that no matter what memory length we apply, either we see a nearly unbiased behavior or the first occurrence visibly attracts further repetitions. 
Our Theorem~1 asserts that if we observe longer and longer blocks $B$, repelling from any distance must decay in both 
measure and intensity (while attracting can persist), so that for majority of long blocks we will see the behavior as 
described above.

We also note, that by pushing the graph of $\tilde F_B$ downward (compared to $1-e^{-t}$), attracting contributes 
to increasing the expected value of the associated random variable, i.e., of the hitting time. In case of 
attracting assisted by only very small intensity repelling, the average waiting time for the first occurrence 
of the event $B$ is increased in comparison to unbiased (may even not exist). Thus, instinctively judging the 
probability of the event by (the inverse of) the waiting time for the first occurrence we will typically 
underestimate it. All the more we are surprised, when the following occurrences happen after a considerably 
shorter time. This additionally strengthens the phenomenon's appearance.

Another consequence of attracting not assisted by repelling (or assisted by repelling of a very small intensity) 
is an increased variance of the return time statistic (the variance may even cease to exist). Thus, again, 
the gaps between the occurrences of $B$ are driven away from the expected value, toward the extremities 0
and $\infty$, and hence, into the pattern of clusters separated by longer pauses. We skip the elementary 
estimations of the variance. 

It must be reminded: we do not claim, for any class of processes, that occurrences of long blocks will actually 
deviate from unbiased. There are conditions, weaker than full independence, under which the distributions of 
the normalized return times of long blocks converge almost surely to the exponential law. It is so, for instance,
in Markov processes (with finite memory). In fact, such convergence is implied by a sufficient rate of mixing 
([A-G], [H-S-V]). Yet, such processes seem to be somewhat exceptional and we expect that attracting rules in
majority of processes (see the Question 5 at the end of the paper). As we have already mentioned, at least that much is true, that in any dynamical system with positive entropy there exist partitions with strong attracting properties.

\medskip
It is important not to be misled by an oversimplified approach. The ``decay of repelling'' 
in positive entropy processes appears to agree with the intuitive understanding of entropy as chaos: 
repelling is a ``self-organizing'' property; it leads to a more uniform, hence less chaotic, distribution 
of an event along a typical orbit. Thus one might expect that repelling with intensity $\epsilon$ revealed 
by a fraction $\xi$ of all $n$-blocks contributes to lowering an upper estimate of the entropy by some 
percentage proportional to $\xi$ and depending increasingly on $\epsilon$. If this happens for infinitely 
many lengths $n$ with the same parameters $\xi$ and $\epsilon$, the entropy should be driven to zero by a 
geometric progression. Surprisingly, it is not quite so, 
and the phenomenon has more subtle grounds. We will present an example which exhibits the incorrectness of 
such intuition (see the Example~1 and the preceding discussion). Also, it will become obvious from the proof, 
that there is no gradual reduction of the entropy. The entropy is ``killed completely in one step'', that 
means, positive entropy and persistent repelling lead to a contradiction by examining the blocks of one 
sufficiently large length $n$; we do not use any iterated procedure requiring repelling for infinitely many 
lengths.

\heading Notation and preliminary facts\endheading
We now establish further notation and preliminaries needed in the proof. If $\Bbb A\subset\z$ then we
will write $\Cal P^{\Bbb A}$ to denote the partition or sigma-field $\bigvee_{i\in \Bbb A}\sigma^{-i}(\Cal P)$.
We will abbreviate $\Cal P^n = \Cal P^{[0,n)}$, $\Cal P^{-n} = \Cal P^{[-n,-1]}$, $\Cal P^- = \Cal P^{(-\infty,-1]}$ 
(a ``finite future'', a ``finite past'', and the ``full past'' of the process).

We assume familiarity of the reader with the basics of entropy for finite partitions and sigma-fields
in a standard probability space. Our notation is compatible with [P] and we refer the reader to 
this book, as well as [Sh] and [W], for background and proofs. In particular, we will be using the following:
\item * The entropy of a partition equals $H(\Cal P) = -\sum_{A\in\Cal P}\mu(A)\log_2(\mu(A))$.
\item * For two finite partitions $\Cal P$ and $\Cal B$, the conditional entropy $H(\Cal P|\Cal B)$ 
is equal to $\sum_{B\in\Cal B}\mu(B)H_B(\Cal P)$, where $H_B$ is the entropy evaluated for the conditional 
measure $\mu_B$ on $B$.
\item * The same formula holds for conditional entropy given a sub-sigma-field $\Cal C$, i.e., 
$$
\sum_{B\in\Cal B}\mu(B)H_B(\Cal P|\Cal C)=H(\Cal P|\Cal B\vee\Cal C). 
$$ 
\smallskip
\item * The entropy of the process is given by any one of the formulas below
$$
h = H(\Cal P|\Cal P^-) = \tfrac 1r H(\Cal P^r|\Cal P^-) = \lim_{r\to\infty}\tfrac 1r H(\Cal P^r).
$$

We will exploit the notion of $\epsilon$-independence for partitions and sigma-fields. The definition below
is an adaptation from [Sh], where it concerns finite partitions only. See also [Sm] for treatment of 
countable partitions. Because ``$\epsilon$'' is reserved for the intensity of repelling, we will speak 
about $\beta$-independence.

\definition{Definition 2} Fix $\beta>0$. 
A partition $\Cal P$ is said to be {\it $\beta$@-independent} of a sigma-field $\Cal B$ if for any
$\Cal B$-measurable countable partition $\Cal B'$ holds 
$$
\sum_{A\in\Cal P, B\in\Cal B'}|\mu(A\cap B)-\mu(A)\mu(B)|\le\beta.
$$ 
A process \proc\ is called an {\it $\beta$-independent process} if $\Cal P$ is $\beta$-independent
of the past $\Cal P^-$.
\enddefinition

A partition $\Cal P$ is independent of another partition or a sigma-field $\Cal B$ if and only if 
$H(\Cal P|\Cal B) = H(\Cal P)$. The following approximate version of this fact holds (see [Sh, Lemma 7.3]
for finite partitions, from which the case of a sigma-field is easily derived). 

\proclaim{Fact 1}
A partition $\Cal P$ is $\beta$-independent of another partition or a sigma-field $\Cal B $ if 
$H(\Cal P|\Cal B)\ge H(\Cal P)-\xi$, for $\xi$ sufficiently small. \qed
\endproclaim
In course of the proof, a certain lengthy condition will be in frequent use. Let us introduce an abbreviation:

\definition{Definition 3}
Given a partition $\Cal P$ of a space with a probability measure $\mu$ and $\delta>0$, we will say that a 
property $\Phi(A)$ {\it holds for $A\in\Cal P$ with $\mu$-tolerance $\delta$} if 
$$
\mu\left(\bigcup\{A\in\Cal P: \Phi(A)\}\right)\ge 1-\alpha.
$$
\enddefinition

We shall also need an elementary estimate, whose proof is an easy exercise.

\proclaim{Fact 2}
For each $A\in\Cal P$, $H(\Cal P) \le (1-\mu(A))\log_2(\#\Cal P) + 1$. \qed
\endproclaim

In addition to the random variables of absolute and normalized return times $R_B$ and $\overline R_B$,
we will also use the analogous notions of the $k^{\text{th}}$ absolute return time
$$
R^{(k)}_B = \min\{i: \#\{0<j\le i: \sigma^j(y)\in B\} = k\},
$$ 
and of the normalized $k^{\text{th}}$ return time $\overline R^{(k)}_B = \mu(B)R^{(k)}_B$ (both defined on $B$), 
with $F^{(k)}_B$ always denoting the distribution function of the latter. Clearly, the expected value of 
$\overline R^{(k)}_B$ equals $k$. 

\heading The idea of the proof and the basic lemma\endheading
Before we pass to the formal proof of Theorem 1, we would like to have the reader oriented in the mainframe 
of the idea behind it. We intend to estimate (from above, by $1-e^{-t}+\epsilon$) the function $G_{B\!A}$, 
for long blocks of the form $B\!A\in\Cal P^{[-n,r)}$.
The ``positive'' part $A$ has a fixed length $r$, while we allow the ``negative'' part $B$ to be arbitrarily long.
There are two key ingredients leading to the estimation. The first one, contained in Lemma 3, is the observation 
that for a fixed typical $B\in \Cal P^{-n}$, the part of the process induced on $B$ (with the conditional 
measure $\mu_B$) generated by the partition $\Cal P^r$, is not only a $\beta$-independent process, but it is 
also $\beta$-independent of many returns times $R_B^{(k)}$ of the cylinder $B$ (see the Figure 2).    

$$
\align
^{\text{coordinate }0}&_\downarrow \\ \vspace{-5pt}
...\boxed{\ \ _{_{B\phantom{_1}}} \ \ }\boxed{_{_{A_{\text{-}1}}}}...............\boxed{\ \ _{_{B\phantom{_1}}} \ \ }
&\boxed{_{_{A_0}}}..\boxed{\ \ _{_{B\phantom{_1}}} \ \ }\boxed{_{_{A_1}}}..........\boxed{\ \ _{_{B\phantom{_1}}} \ \ }\boxed{_{_{A_2}}}....\boxed{\ \
_{_{B\phantom{_1}}} \ \ }\boxed{_{_{A_3}}}....
\endalign
$$
{\eightpoint\it Figure 2: The process $\dots A_{-1}A_0A_1A_2\dots$ of $r$-blocks following the copies of 
$B$ is a $\beta$@-independent process with additional $\beta$-independence properties from the positioning 
of the copies of $B$.}

\medskip\noindent
This allows us to decompose (with high accuracy) the distribution function $F_{B\!A}$ of the normalized return 
time of $B\!A$ as follows:
$$
\gather
F_{B\!A}(t) = \mu_{B\!A}\{\overline R_{B\!A} \le t\} = \mu_{B\!A}\{R_{B\!A} \le \tfrac t{\mu(B\!A)}\} = \\
\sum_{k\ge 1} \mu_{B\!A}\{R^{(B)}_A = k, R^{(k)}_B \le \tfrac t{p\mu(B)}\}\approx 
\sum_{k\ge 1}\mu_{B\!A}\{R^{(B)}_A=k\}\cdot\mu_B\{\overline R^{(k)}_B \le \tfrac tp\} 
\approx \\ \sum_{k\ge 1}p(1-p)^{k-1}\cdot F^{(k)}_B(\tfrac tp),
\endgather
$$
where $R^{(B)}_A$ denotes the first (absolute) return time of $A$ in the process induced on $B$, and $p = \mu_B(A)$.

The second key observation is, assuming for simplicity full independence, that when trying to model some 
repelling for the blocks $B\!A$, we ascertain that it is largest, when the occurrences of $B$ are purely 
periodic. Any deviation from periodicity of the $B$'s may only lead to increasing the intensity of attracting 
between the copies of $B\!A$, never that of repelling. We will explain this phenomenon more formally in 
a moment. Now, if $B$ does appear periodically, then the normalized return time of $B\!A$ is governed by 
the same geometric distribution as the normalized return time of $A$ in the independent process induced
on $B$. If $p$ is small, this geometric distribution function becomes nearly the unbiased exponential 
law $1-e^{-t}$. The smallness of $p$ is {\it a priori} regulated by the choice of the parameter $r$ (Lemma~1). 

The phenomena that, assuming full independence, the repelling of $B\!A$ is maximized by periodic 
occurrences of $B$, and that even then there is nearly no repelling, is captured by the following 
elementary lemma, which will be also useful later, near the end of the rigorous proof.

\proclaim{Lemma 0}
Fix some $p\in(0,1)$. Let $F^{(k)}$ ($k\ge 1$) be a \sq\ of distribution functions on 
$[0,\infty)$ such that the expected value of the distribution associated to $F^{(k)}$ equals $k$. Define   
$$
F(t) = \sum_{k\ge 1}p(1-p)^{k-1} F^{(k)}(\tfrac tp),
\text{ \ \ and \ \ }
G(t) = \int_0^t 1-F(s)ds.
$$
Then $G(t)\le \frac1{\log e_p}(1-e_p^{-t})$, where $e_p = (1-p)^{-\frac 1p}$.
\endproclaim

\demo{Proof} We have
$$
G(t) = \sum_{k\ge 1}p(1-p)^{k-1}\int_0^t 1- F^{(k)}(\tfrac sp) ds.
$$
We know that $F^{(k)}(t)\in [0,1]$ and that $\int_0^\infty 1- F^{(k)}(s) ds = k$ (the expected value). 
With such constraints, it is the indicator function $1_{[k,\infty)}$ that maximizes the integrals from 
$0$ to $t$ simultaneously for every $t$ (because the ``mass'' $k$ above the graph is for such choice 
of the function swept maximally to the left). The rest follows by direct calculations:
$$
\gather
G(t) \le \sum_{k\ge 1}p(1-p)^{k-1}\int_0^t 1_{[0,k)}(\tfrac sp) ds = 
\int_0^t \sum_{k=\lceil\frac sp\rceil}^\infty p(1-p)^{k-1} ds =\\ \int_0^t (1-p)^{\lceil\frac sp\rceil} ds \le \frac{(1-p)^{\frac tp}-1}{\log (1-p)^{\frac 1p}}.\qed
\endgather
$$ 
\enddemo

Recall that the maximizing distribution functions $F_B^{(k)} = 1_{[k,\infty)}$ occur, for the normalized 
return time of a set $B$, precisely when $B$ is visited periodically. This explains our former statement on
this subject.

\medskip
Let us comment a bit more on the first key ingredient, the $\beta$-independence. Establishing it is the most 
complicated part of the argument. The idea is to prove conditional (given a ``finite past'' $\Cal P^{-n}$) 
$\beta$-independence of the ``present'' $\Cal P^r$ from jointly the full past and a large part of the future,
responsible for the return times of majority of the blocks $B\in\Cal P^{-n}$. But the future part must not be 
too large. Let us mention the existence of ``bilaterally deterministic'' processes with positive entropy 
(first discovered by Gurevi\v c [G], see also [O-W1]), in which the sigma-fields generated by the coordinates 
$(-\infty,-m]\cup[m,\infty)$ do not decrease with $m$ to the Pinsker factor; they are all equal to the entire 
sigma-field. (Coincidently, our Example~1 has precisely this property; see the Remark 2.) Thus, in order 
to maintain any trace of independence of the ``present'' from our sigma-field
already containing the entire past, its part in the future must be selected with an extreme care. Let us also 
remark that an attempt to save on the future sigma-fields by adjusting them 
individually to each block $B_0\in\Cal P^{-n}$ falls short, mainly because of the ``off diagonal effect''; 
suppose $\Cal P^r$ is conditionally (given $\Cal P^{-n}$) nearly independent of a sigma-field which determines 
the return times of only one selected block $B_0\in \Cal P^{-n}$. The independence still holds conditionally 
given any cylinder $B\in\Cal P^{-n}$ from a collection of a large measure, but unfortunately, this collection 
can always miss the selected cylinder $B_0$. In Lemmas 2 and 3, we succeed in finding a sigma-field (containing 
the full past and a part of the future), of which $\Cal P^r$ is conditionally $\beta$-independent, and which  
``nearly determines'', for majority of blocks $B\in\Cal P^{-n}$, some finite number of their sequential return 
times (probably not all of them). This finite number is sufficient to allow the described earlier 
decomposition of the distribution function $F_{B\!A}$. 

\heading The proof \endheading
Throughout the sequel we assume ergodicity and that the entropy $h$ of \proc\ is 
positive. We begin our computations with an auxiliary lemma allowing us to assume 
(by replacing $\Cal P$ by some $\Cal P^r$) that the elements of the ``present'' partition are small, 
relatively in most of $B\in\Cal P^n$ and for every $n$. Note that the Shannon-McMillan-Breiman Theorem
is insufficient: for the conditional measure the error term depends increasingly on $n$, which we do not fix.

\proclaim{Lemma 1} For each $\delta$ there exists an $r\in\na$ such that 
for every $n\in\na$ the following holds for $B\in\Cal P^{-n}$ with $\mu$-tolerance $\delta$:
$$
\text{for every }A\in\Cal P^r,\ \ \mu_B(A)\le \delta.
$$
\endproclaim

\demo{Proof} 
Let $\alpha$ be so small that 
$$
\sqrt\alpha\le\delta \text{\ \ and\ \ } \frac{h-3\sqrt\alpha}{h+\alpha}\ge 1-\frac\delta2,
$$
and set $\gamma = \frac\alpha{\log_2(\#\Cal P)}$. Let $r$ be so big that 
$$
\frac 1r \le \alpha, \ \ \frac1{r(h+\alpha)}\le\frac\delta2,
$$
and that there exists a collection $\overline{\Cal P^r}$ of no more than $2^{r(h+\alpha)}-1$ elements of 
$\Cal P^r$ whose joint measure $\mu$ exceeds $1-\gamma$ (by the Shannon-McMillan-Breiman Theorem).

Let $\widetilde{\Cal P^r}$ denote the partition into the elements of $\overline{\Cal P^r}$ and the complement 
of their union, and let $\Cal R$ be the partition into the remaining elements of $\Cal P^r$ and the complement 
of their union, so that $\Cal P^r = \widetilde{\Cal P^r}\vee \Cal R$. For any $n$ we have 
$$
\gather
rh = H(\Cal P^r|\Cal P^-)\le H(\Cal P^r|\Cal P^{-n}) = 
H(\widetilde{\Cal P^r}\vee \Cal R|\Cal P^{-n})=\\
H(\widetilde{\Cal P^r}|\Cal R \vee \Cal P^{-n})+ H(\Cal R|\Cal P^{-n})\le
H(\widetilde{\Cal P^r}|\Cal P^{-n}) + H(\Cal R) \le  \\
\sum_{B\in\Cal P^{-n}}\mu(B)H_B(\widetilde{\Cal P^r}) + \gamma r\log_2(\#\Cal P) + 1
\endgather
$$
(we have used Fact 2 for the last passage). After dividing by $r$, we obtain
$$
\sum_{B\in\Cal P^{-n}}\mu(B)\tfrac 1r H_B(\widetilde{\Cal P^r}) \ge h - \gamma\log_2(\#\Cal P) - \tfrac1r 
\ge h - 2\alpha.
$$
Because each term $\tfrac 1r H_B(\widetilde{\Cal P^r})$ is not larger than $\frac 1r\log_2(\#\widetilde{\Cal P^r})$ 
which was set to be at most $h+\alpha$, we deduce that 
$$ 
\tfrac1r H_B(\widetilde{\Cal P^r}) \ge h - 3\sqrt\alpha
$$
holds for $B\in\Cal P^{-n}$ with $\mu$-tolerance $\sqrt\alpha$, hence also with $\mu$-tolerance $\delta$. 
On the other hand, by Fact 2, for any $B$ and $A\in\widetilde{\Cal P^r}$, holds: 
$$
H_B(\widetilde{\Cal P^r}) \le (1-\mu_B(A))\log_2(\#\widetilde{\Cal P^r})+1 \le (1-\mu_B(A))r(h+\alpha)+1.
$$
Combining the last two displayed inequalities we establish that, with $\mu$-tolerance $\delta$ for $B\in\Cal P^{-n}$ 
and then for every $A\in\widetilde{\Cal P^r}$, holds
$$
1-\mu_B(A)\ge \frac{h-3\sqrt\alpha}{h+\alpha} -\frac1{r(h+\alpha)}\ge 1-\delta.
$$
So, $\mu_B(A)\le \delta$. Because $\Cal P^r$ refines $\widetilde{\Cal P^r}$, the elements of $\Cal P^r$ 
are not larger.
\qed\enddemo

We continue the proof with a lemma which can be deduced from [R1,~Lemma~3]. 
We provide a direct proof. For $\alpha>0$ and $M\in\na$ let 
$$
S(M,\alpha) = \bigcup_{m\in\z}[mM+\alpha M, (m+1)M-\alpha M)\cap\z.
$$

\proclaim{Lemma 2}
For fixed $\alpha$ and $r$ there exists $M_0$ such that for every $M\ge M_0$ holds, 
$$
H(\Cal P^r|\Cal P^-\vee \Cal P^{S(M,\alpha)})\ge rh-\alpha
$$ 
(see the Figure 3).
$$
***********\circ\circ..************..........************..........************..........
$$
\endproclaim\noindent
{\eightpoint\it Figure 3. The circles indicate the coordinates 0 through $r-1$, the conditioning sigma-filed 
is over the coordinates marked by stars, which includes the entire past and part of the future with gaps of size
$2\alpha M$ repeated periodically with period $M$ (the first gap is half the size).}

\demo{Proof} First assume that $r=1$. 
Denote also
$$
S'(M,\alpha) = \bigcup_{m\in\z}[mM+\alpha M, (m+1)M)\cap\z.
$$
Let $M$ be so large that $H(\Cal P^{(1-\alpha)M})< (1-\alpha)M(h+\gamma)$, where 
$\gamma = \frac{\alpha^2}{2(1-\alpha)}$. Then, for any $m\ge 1$,
$$
H(\Cal P^{S'(M,\alpha)\cap[0,mM)}|\Cal P^-)\le 
H(\Cal P^{S'(M,\alpha)\cap[0,mM)})<(1-\alpha)mM(h+\gamma).
$$ 
Because $H(\Cal P^{[0,mM)}|\Cal P^-) = mMh$, the complementary part of entropy must exceed 
$mMh - (1-\alpha)mM(h+\gamma)$ (which equals $\alpha mM (h-\tfrac\alpha2)$), i.e., we have
$$
H(\Cal P^{[0,mM)\setminus S'(M,\alpha)}|\Cal P^-\vee \Cal P^{S'(M,\alpha)\cap[0,mM)})
> \alpha mM (h-\tfrac\alpha2).
$$
Breaking the last entropy term as a sum over $j\in [0,mM)\setminus S'(M,\alpha)$ of the conditional
entropies of $\sigma^{-j}(\Cal P)$ given the sigma-field over all coordinates left of $j$ and all coordinates
from $S'(M,\alpha)\cap[0,mM)$ right of $j$, and because every such term is at most $h$, we deduce 
that more than half of these terms reach or exceed $h-\alpha$. So, a term not smaller than $h-\alpha$ 
occurs for a $j$ within one of the gaps in the left half of $[0,mM)$.
Shifting by $j$, we obtain $H(\Cal P|\Cal P^-\vee \sigma^i(\Cal P^{S'(M,\alpha)\cap [0,\frac{mM}2)}))\ge h-\alpha$,
where $i\in [0,\alpha M)$ denotes the relative position of $j$ in the gap. As we increase $m$, one value 
$i$ will repeat in this role along a subsequence $m'$. The operation $\vee$ is continuous for 
increasing \sq s of sigma-fields, hence $\Cal P^-\vee \sigma^i(\Cal P^{S'(M,\alpha)\cap [0,\frac{m'M}2)})$ 
converges over $m'$ to $\Cal P^-\vee \sigma^i(\Cal P^{S'(M,\alpha)})$. The entropy is continuous for such 
passage, hence $H(\Cal P|\Cal P^-\vee \sigma^i(P^{S'(M,\alpha)})\ge h-\alpha$. The assertion now follows 
because $S(M,\alpha)$ is contained in $S'(M,\alpha)$ shifted to the left by any $i\in[0,\alpha M)$.

Finally, if $r>1$, we can simply argue for $\Cal P^r$ replacing $\Cal P$. This will impose that
$M_0$ and $M$ are divisible by $r$, but it is not hard to see that for large $M$ the argument
works without divisibility at a cost of a slight adjustment of $\alpha$.
\qed\enddemo

For a long block $B\in\Cal P^{-n}$ let \procbr\ denote the process induced on $B$ generated by the restriction 
$\Cal P^r_B$ of $\Cal P^r$ to $B$ ($\sigma_B$ is the first return time map on $B$). 
The following lemma is the crucial item in our argument.

\proclaim{Lemma 3}
For every $\beta>0$, $r\in\na$ and $K\in\na$ there exists $n_0$ such that for every $n\ge n_0$, 
with $\mu$-tolerance $\beta$ for $B\in\Cal P^{-n}$, with respect to $\mu_B$, $\Cal P^r$ is $\beta$-independent 
of jointly the past $\Cal P^-$ and the first $K$ return times to $B$, $R^{(k)}_B$ ($k\in[1,K]$). 
In particular, \procbr\ is a $\beta$-independent process.
\endproclaim

\demo{Proof} 
We choose $\xi$ according to Fact 1, so that $\frac\beta2$-independence is implied.  Let $\alpha$ satisfy 
$$
0<\tfrac{2\alpha}{h-\alpha}<1, \ \ 18K\sqrt\alpha<1,\ \ \sqrt{2\alpha}<\xi,\ \ K\root 4\of\alpha<\tfrac\beta2.
$$ 
Let $n_0$ be so large that $H(\Cal P^r|\Cal P^{-n})<rh+\alpha$ for every $n\ge n_0$ and that for every 
$k\in [1,K]$ with $\mu$-tolerance $\alpha$ for $B\in\Cal P^{-n}$ holds
$$
\mu_B\{2^{n(h-\alpha)}\le R^{(k)}_B\le 2^{n(h+\alpha)}\}>1-\alpha
$$
(we are using Ornstein-Weiss Theorem [O-W2]; the multiplication by $k$ is consumed by $\alpha$ in the exponent). 
Let $M_0\ge 2^{n_0(h-\alpha)}$ be so large that the assertion of Lemma~2 holds for $\alpha$, $r$ and $M_0$, 
and that for every $M\ge M_0$, 
$$
(M+1)^{1+\frac{2\alpha}{h-\alpha}}<\alpha M^2 \text{ \ and \ }\tfrac{\log_2(M+1)}{M(h-\alpha)}<\alpha. 
$$
We can now redefine (enlarge) $n_0$ and $M_0$ so that $M_0=\lfloor 2^{n_0(h-\alpha)}\rfloor$. Similarly,
for each $n\ge n_0$ we set $M_n=\lfloor 2^{n(h-\alpha)}\rfloor$. Observe, that the interval where the first 
$K$ returns of most $n$-blocks $B$ may occur (up to probability $\alpha$), is contained in $[M_n, \alpha M_n^2]$
(because $2^{n(h+\alpha)}\ \le (M_n+1)^{1+\frac{2\alpha}{h-\alpha}}<\alpha M_n^2$). 

At this point we fix some 
$n\ge n_0$. The idea is to carefully select an $M$ between $M_n$ and $2M_n$ (hence not smaller than $M_0$), 
such that the initial $K$ returns of nearly every $n$-block happen most likely inside (with all its $n$ 
symbols) the set $S(M,\alpha)$, so that they are ``controlled'' by the sigma-field $\Cal P^{S(M,\alpha)}$. 
Let $\alpha' = \alpha + \frac n{M_n}$, so that every $n$-block overlapping with $S(M,\alpha')$ is completely
covered by $S(M,\alpha)$. By the second assumption on $M\ge M_0$ and by the formula connecting $M_n$ and $n$, 
we have $\alpha'<2\alpha$. To define $M$ we will invoke the triple Fubini Theorem. Fix $k\in[1,K]$ and consider 
the probability space 
$$
\Cal P^{-n}\times [M_n,2M_n]\times \na
$$ 
equipped with the (discrete) measure $\Cal M$ whose marginal on $\Cal P^{-n}\times [M_n,2M_n]$ is the product 
of $\mu$ (more precisely, of its projection onto $\Cal P^{-n}$) with the uniform distribution on the integers 
in $[M_n,2M_n]$, while, for fixed $B$ and $M$, the measure on the corresponding $\na$-section is the distribution 
of the random variable $R_B^{(k)}$. In this space let $S$ be the set whose $\na$-section for a fixed 
$M$ (and any fixed $B$) is the set $S(M,\alpha')$. We claim that for every $l\in [M_n,\alpha M_n^2]\cap\na$ 
(and any fixed $B$) the $[M_n,2M_n]$-section of $S$ has measure exceeding $1-16\alpha$. This is quite obvious 
(even for every $l\in [M_n,\infty)$ and with $1-15\alpha$) if $[M_n,2M_n]$ is equipped with the normalized 
Lebesgue measure (see the Figure 4). 

\input epsf.tex
\epsfxsize=13truecm
\centerline{\epsfbox{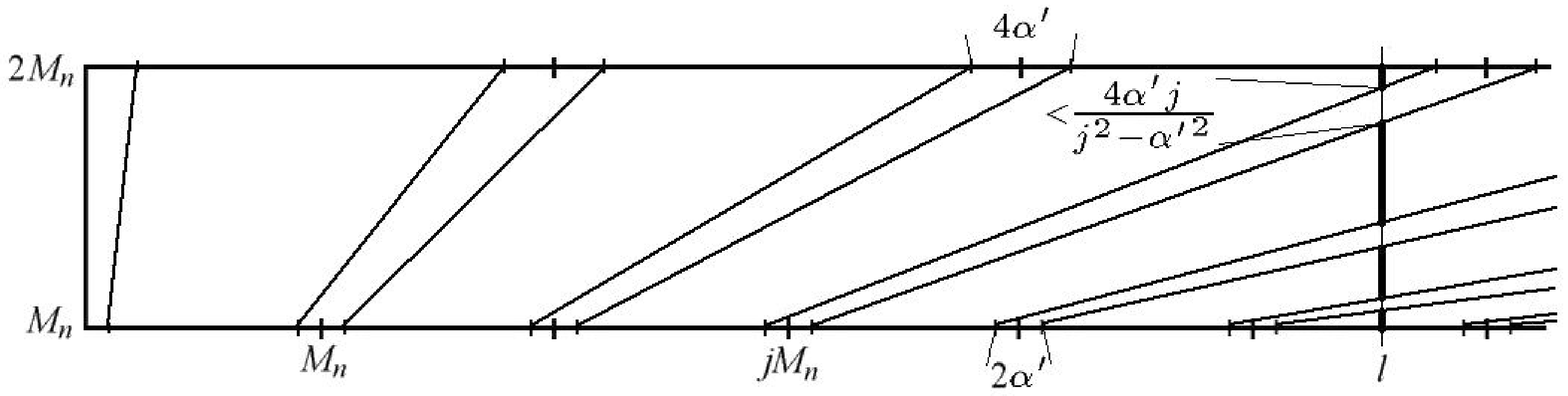}}
{\eightpoint\it Figure 4: The complement of $S$ splits into thin skew strips shown in the picture.
The normalized Lebesgue measure of any vertical section of the $j^{\text{\rm th}}$ strip (starting at $jM_n$ 
with $j\ge 1$) is at most $\frac{4\alpha'j}{j^2 - {\alpha'}^2}\le \frac{5\alpha'}j\le \frac{10\alpha}j$. 
Each vertical line at $l\ge M_n$ intersects strips with indices $j, j+1, j+2$ up to at most $2j$ 
(for some $j$), so the joint measure of the complement of the section of $S$ does not exceed $15\alpha$.}

\medskip
$$
\align \ \ 
&\ \ \ \ \ \ \ \ \ \ \ \ \ \ \ \ \ \ S\\
&\ \ \ \ \ \ \ \ \ \ \ \ \swarrow\ \ \downarrow\ \ \ \searrow\\
\vspace{-5pt}
&|\!\overline{..}.\overline{...}.\overline{....}.\overline{.....}..\overline{......}..
\overline{.......}..\overline{........}...\overline{.........}...\overline{..........}...\overline{..........}...\overline{.....}\!| \\
&\!M_n\hphantom{...................................................................................}2M_n
\endalign
$$
{\eightpoint\it Figure 5: The discretization replaces the Lebesgue measure by the uniform measure on 
$M_n$ integers, thus the measure of any interval can deviate from its Lebesgue measure by at most $\frac 1{M_n}$.
For $l\le \alpha M_n^2$ the corresponding section of $S$ (in this picture drawn horizontally) 
consists of at most $\alpha M_n$ intervals, so its measure can deviate by no more than $\alpha$.}

In the discrete case, however, {\it a priori} it might happen that the integers
along some $[M_n,2M_n]$-section often ``miss'' the section of $S$ leading to a decreased measure value. (For 
example, it is easy to see that for $l = (2M_n)!$ the measure of the section of $S$ is zero.) But because we 
restrict to $l\le \alpha M_n^2$, the discretization does not affect the measure of the section of $S$ 
by more than $\alpha$, and the estimate with $1-16\alpha$ holds (see the Figure 5 above).

Taking into account all other inaccuracies (the smaller than $\alpha$ part of $S$ outside $[M_n,\alpha M_n^2]$ 
and the smaller than $\alpha$ part of $S$ projecting onto blocks $B$ which do not obey the Ornstein-Weiss return 
time estimate) it is safe to claim that 
$$
\Cal M(S)>1-18\alpha.
$$ 
This implies that for every $M$ from a set of measure at least $1-18\sqrt{\alpha}$ 
the measure of the $(\Cal P^{-n}\times \na)$-section of $S$ is larger than or equal to $1-\sqrt\alpha$. 
For every such $M$, with $\mu$-tolerance $\root4\of\alpha$ for $B\in\Cal P^{-n}$, the probability $\mu_B$ 
that the $k^{\text{\rm th}}$ repetition of $B$ falls in $S(M,\alpha')$ (hence with all its $n$ terms inside 
the set $S(M,\alpha)$) is at least $1-\root 4\of\alpha$.

Because $18K\sqrt\alpha<1$, there exists at least one $M$ for which the above holds for every $k\in[1,K]$.
This is our final choice of $M$. For this $M$, with $\mu$-tolerance $K\root 4\of\alpha$, all considered $K$ 
returns of $B$ are, with probability $1-\root 4\of\alpha$ (each), determined by the sigma-field 
$\Cal P^{S(M,\alpha)}$. More precisely, for each $k\in[1,K]$ there is a set $U_k$ of measure $\mu_B$ at most 
$\root4\of\alpha$ such that the sets $\{R^{(k)}_B = i\}\setminus U_k$ agree with some $Q^{(k)}_i\setminus U_k$, 
where each $Q^{(k)}_i$ is $\Cal P^{S(M,\alpha)}$-measurable. Thus, we can modify the variable $R^{(k)}_B$ so it is 
$\Cal P^{S(M,\alpha)}$-measurable and equal to the original except on $U_k$.
We denote such a modification by $\tilde R^{(k)}_B$.

Let us go back to our entropy estimates. We have, by Lemma 2,
$$
\gather
\sum_{B\in\Cal P^{-n}}\mu(B)H_B(\Cal P^r|\Cal P^-\vee \Cal P^{S(M,\alpha)}) = 
H(\Cal P^r|\Cal P^{-n}\vee\Cal P^-\vee \Cal P^{S(M,\alpha)})=\\
H(\Cal P^r|\Cal P^-\vee \Cal P^{S(M,\alpha)}) \ge
rh -\alpha \ge H(\Cal P^r|\Cal P^{-n})-2\alpha =\\
\sum_{B\in\Cal P^{-n}}\mu(B)H_B(\Cal P^r)-2\alpha.
\endgather
$$
Because $H_B(\Cal P^r|\Cal P^-\vee \Cal P^{S(M,\alpha)})\le H_B(\Cal P^r)$ for every $B$,
we deduce that with $\mu$-tolerance $\sqrt{2\alpha}$ for $B\in\Cal P^{-n}$ must hold 
$$
H_B(\Cal P^r|\Cal P^-\vee \Cal P^{S(M,\alpha)})\ge H_B(\Cal P^r)-\sqrt{2\alpha}\ge H_B(\Cal P^r)-\xi.
$$
Combining this with the preceding arguments, with $\mu$-tolerance $K\root 4\of\alpha + \sqrt{2\alpha}<\beta$
for $B\in\Cal P^{-n}$ both the above entropy inequality holds, and we have the $\Cal P^{S(M,\alpha)}$@-mea\-surable 
modifications $\tilde R^{(k)}_B$ of the return times. By the choice of $\xi$, we obtain that with respect to $\mu_B$, $\Cal P^r$ is jointly $\frac\beta2$-independent of the past and the modified return times 
$\tilde R^{(k)}_B$ ($k\in [1,K]$). Because $\mu(\bigcup_{k\in[1,K]} U_k)\le K\root4\of\alpha<\frac\beta2$, 
this clearly implies $\beta$-independence if each $\tilde R^{(k)}_B$ is replaced by $R^{(k)}_B$.
\qed\enddemo

To complete the proof of Theorem 1 it now remains to put the items together.

\demo{Proof of Theorem 1}
Fix an $\epsilon>0$. On $[0,\infty)$, the functions 
$$
g_p(t) = \min\{1,\tfrac1{\log e_p}(1-e_p^{-t}) +pt\},
$$ 
where $e_p = (1-p)^{-\frac 1p}$, decrease uniformly to $1-e^{-t}$ as $p\to 0^+$. So, let $\delta$ be such 
that $g_{\delta}(t)\le 1-e^{-t}+\epsilon$ for every $t$. We also assume that 
$$
(1-2\delta)(1-\delta)\ge 1-\epsilon.
$$
Let $r$ be specified by Lemma 1, so that $\mu_B(A)\le \delta$ for every $n\ge 1$, every $A\in\Cal P^r$
and for $B\in\Cal P^{-n}$ with $\mu$-tolerance $\delta$. On the other hand, once $r$ is fixed, the 
partition $\Cal P^r$ has at most 
$(\#\Cal P)^r$ elements, so with $\mu_B$-tolerance $\delta$ for $A\in\Cal P^r$, $\mu_B(A)\ge \delta(\#\Cal P)^{-r}$. 
Let $\Cal A_B$ be the subfamily of $\Cal P^r$ (depending on $B$) where this inequality holds.
Let $K$ be so large that for any $p\ge \delta(\#\Cal P)^{-r}$,
$$
\sum_{k=K+1}^\infty p(1-p)^k <\tfrac\delta2,
$$
and choose $\beta<\delta$ so small that 
$$
(K^2+K+1)\beta<\tfrac\delta2.
$$ 
The application of Lemma 3 now provides an $n_0$ such that for any $n\ge n_0$, with $\mu$-tolerance $\beta$
for $B\in\Cal P^{-n}$, the process induced on $B$ generated by $\Cal P^r$ has the desired $\beta$-independence 
properties involving the initial $K$ return times of $B$. So, with tolerance $\delta+\beta<2\delta$ we have 
both, the above $\beta$-independence and the estimate $\mu_B(A)<\delta$ for every $A\in\Cal P^r$. Let $\Cal B_n$ 
be the subfamily of $\Cal P^{-n}$ where these two conditions hold. Fix some $n\ge n_0$.

Let us consider a cylinder set $B\cap A\in\Cal P^{[-n,r)}$ (or, equivalently, the block $B\!A$), 
where $B\in\Cal B_n$, $A\in \Cal A_B$. The length of $B\!A$ is $n+r$, which represents an arbitrary integer 
larger than $n_0+r$. Notice that the family of such sets $B\!A$ covers more than $(1-2\delta)(1-\delta)\ge 1-\epsilon$ 
of the space. 

We will examine the distribution of the normalized first return time for $B\!A$. 
In addition to our customary notations of return times, let $R^{(B)}_A$ be the first (absolute) 
return time of $A$ in \procbr, i.e., the variable defined on $B\!A$, counting the number of visits to 
$B$ until the first return to $B\!A$. Let $p = \mu_B(A)$. We have
$$
\gather
F_{B\!A}(t) = \mu_{B\!A}\{\overline R_{B\!A} \le t\} = \mu_{B\!A}\{R_{B\!A} \le \tfrac t{\mu(B\!A)}\} = \\
\sum_{k\ge 1} \mu_{B\!A}\{R^{(B)}_A = k, R_B^{(k)} \le \tfrac t{p\mu(B)}\}.
\endgather
$$ 
The $k^{\text{th}}$ term of this sum equals
$$
\tfrac 1p \mu_B(\{A_k = A\}\cap\{A_{k-1}\neq A\}\cap\dots
\cap\{A_1\neq A\}\cap\{A_0 = A\}\cap\{R^{(k)}_B \le \tfrac t{p\mu(B)}\}),
$$
where $A_i$ is the $r$-block following the $i^{\text{th}}$ copy of $B$ (the counting starts from $0$ at
the copy of $B$ positioned at $[-n,-1]$).

By Lemma 3, for $k\le K$, in this intersection of sets each term is $\beta$-independent of the intersection right from it. So, proceeding from the left, we can replace the probabilities of the
intersections  by products of probabilities, allowing an error of $\beta$. Note that the last term equals 
$\mu_B\{\overline R^{(k)}_B \le \tfrac tp\} = F^{(k)}_B(\tfrac tp)$. Jointly, the inaccuracy will not exceed 
$(K+1)\beta$:
$$
\left|\mu_{B\!A}\{R^{(B)}_A = k, R_B^{(k)} \le \tfrac t{p\mu(B)}\} - p(1-p)^{k-1}F^{(k)}_B(\tfrac tp)\right|
\le (K+1)\beta.
$$
Similarly, we also have $\left|\mu_{B\!A}\{R^{(B)}_A=k\}-p(1-p)^{k-1}\right|\le K\beta$, hence the tail of the series 
$\mu_{B\!A}\{R^{(B)}_A=k\}$ above $K$ is smaller than $K^2\beta$ plus the tail of the geometric series $p(1-p)^{k-1}$,
which, by the fact that $p\ge \delta(\#\Cal P)^{-r}$, is smaller than $\tfrac\delta2$. Therefore
$$
F_{B\!A}(t) \approx \sum_{k\ge 1}p(1-p)^{k-1}F^{(k)}_B(\tfrac tp),
$$
up to $(K^2+K+1)\beta +\tfrac\delta2\le\delta$, uniformly for every $t$. By the application of Lemma~0, 
$G_{B\!A}$ satisfies
$$
G_{B\!A}(t) \le \min\{1, \tfrac1{\log e_p}(1-e_p^{-t}) + \delta t\}\le 
g_{\delta}(t) \le 1-e^t+\epsilon
$$
(because $p\le\delta$).
We have proved that for our choice of $\epsilon$ and an arbitrary length $m\ge n_0+r$, with $\mu$-tolerance 
$\epsilon$ for the cylinders $C\in \Cal P^m$, the intensity of repelling between visits to $C$ is at most 
$\epsilon$. This concludes the proof of Theorem 1.
\qed\enddemo

\heading Consequences for limit laws\endheading
The studies of limit laws for return/hitting time statistics are based on the following approach:
For $x\in\Cal P^\z$ define $F_{x,n} = F_B$ (and $\tilde F_{x,n}=\tilde F_B$), where $B$ is the block $x[0,n)$ 
(or the cylinder in $\Cal P^n$ containing $x$). Because for nondecreasing functions $F:[0,\infty)\to[0,1]$, the 
weak convergence coincides with the convergence at continuity points, and it makes the space of such functions 
metric and compact, for every $x$ there exists a well defined collection of limit distributions for $F_{x,n}$ 
(and for $\tilde F_{x,n}$) as $n\to \infty$. They are called {\it limit laws for the return (hitting) 
times at $x$}. Due to the integral relation ($\tilde F_B\approx G_B$) a \sq\ of return time distributions 
converges weakly if and only if the corresponding hitting time distributions converge pointwise (see [H-L-V]), 
so the limit laws for the return times completely determine those for hitting times and {\it vice versa}. 
A limit law is {\it essential} \,if it appears along some sub\sq\ $(n_k)$ for $x$'s in a set of positive measure. 
In particular, the strongest situation occurs when there exists an almost sure limit law along the full \sq\ $(n)$. 
Most of the results concerning the limit laws, obtained so far, can be classified in three major groups:
a) characterizations of possible essential limit laws for specific zero entropy processes (e.g. [D-M], [C-K];
these limit laws are usually atomic for return times or piecewise linear for hitting times), b) finding classes of 
processes with an almost sure exponential limit law along $(n)$ (e.g. [A-G], [H-S-V]), and c) results concerning 
not essential limit laws, limit laws along sets other than cylinders (see [L]; every probabilistic distribution 
with expected value not exceeding 1 can occur in any process as the limit law for such general return times), or 
other very specific topics. As a consequence of our Theorem 1, we obtain, for the first time, a serious 
bound on the possible essential limit laws for the hitting time statistics along cylinders in the 
general class of ergodic positive entropy processes. The statement (1) below is even slightly stronger, 
because we require, for a subsequence, convergence on a positive measure set, but not necessarily to a common limit.

\proclaim{Theorem 2} Assume ergodicity and positive entropy of the process \proc.
\roster
\item If a sub\sq\ $(n_k)$ is such that
$\tilde F_{x,n_k}$ converge pointwise to some limit laws $\tilde F_x$ on a positive measure set 
$A$ of points $x$, then almost surely on $A$, $\tilde F_x(t)\le 1-e^{-t}$ at each $t\ge0$. 
\item If $(n_k)$ grows sufficiently fast, then there is a full measure set, such that for every $x$
in this set holds: $\limsup_k\, \tilde F_{x,n_k}(t)\le 1-e^{-t}$ at each $t\ge0$.
\endroster
\endproclaim

\demo{Proof}
The implication from Theorem 1 to Theorem 2 is obvious and we leave it to the reader. 
For (2) we hint that $(n_k)$ must grow fast enough to ensure summability of the measures of the sets 
where the intensity of repelling persists.
\qed\enddemo

\heading Examples \endheading
The first construction will show that for each $\delta>0$ and $n\in\na$ there exists $N\in\na$ and an ergodic 
process on $N$ symbols with entropy $\log_2N -\delta$, such that the $n$-blocks from a collection of joint 
measure equal to $\frac 1n$ repel with nearly the maximal possible intensity $e^{-1}$. Because $\delta$ 
can be extremely small compared to $\frac 1n$, this construction illustrates, that there is no ``reduction of entropy'' by an amount proportional to the fraction of blocks which reveal strong repelling. 

\example{Example 1} Let $\Cal P$ be an alphabet 
of a large cardinality $N$. Divide $\Cal P$ into two disjoint subsets, one, denoted $\Cal P_0$, of cardinality 
$N_0 = N2^{-\delta}$ and the relatively small (but still very large) rest which we denote by $\{1,2,\dots,r\}$ 
(we will refer to these symbols as ``markers''). For $i = 1,2,\dots,r$, let $\Cal B_i$ be the collection of all 
$n$-blocks whose first $n-1$ symbols belong to $\Cal P_0$ and the terminal symbol is the marker $i$. The cardinality 
of $\Cal B_i$ is $N_0^{n-1}$. Let $\Cal C_i$ be the collection of all blocks of length $nN_0^{n-1}$ obtained 
as concatenations of blocks from $\Cal B_i$ using each of them exactly once. The cardinality of $\Cal C_i$ 
is $(N_0^{n-1})!$. Let $X$ be the subshift whose points are infinite concatenations of blocks from 
$\bigcup_{i=1}^r\Cal C_i$, in which every block belonging to $\Cal C_i$ is followed by a block from 
$\Cal C_{i+1}$ ($1\le i<r$) and every block belonging to $\Cal C_r$ is followed by a block from $\Cal C_1$.
Let $\mu$ be the shift-invariant measure of maximal entropy on $X$. It is immediate to see that the entropy
of $\mu$ is $\frac 1{nN_0^{n-1}}\log_2((N_0^{n-1})!)$, which, for large $N$, nearly equals $\log_2{N_0} = \log_2N -\delta$. Finally observe that the measure of each $B\in\Cal B_i$ equals $\frac 1{nrN_0^{n-1}}$, the joint measure 
of $\bigcup_{i=1}^r\Cal B_i$ is exactly $\frac1n$, and every block $B$ from this family appears in any 
$x\in X$ with gaps ranging between $\tfrac{1-\frac 1r}{\mu(B)}$ and $\tfrac{1+\frac 1r}{\mu(B)}$, 
revealing strong repelling.
\endexample

\remark{Remark 1}
Viewing blocks of length $nrN_0^{n-1}$ starting with a block from $\Cal C_1$ as a new alphabet, and 
repeating the above construction inductively, we can produce an example (with the measure of maximal entropy
on the intersection of systems created in consecutive steps) with entropy $\log_2N - 2\delta$, in which the 
strong repelling will occur with probability $\frac 1{n_k}$ for infinitely many lengths $n_k$.
\endremark

\remark{Remark 2}
The process described in the above remark is (somewhat coincidently; it was not designed for that)
bilaterally deterministic: for every $m\in\na$ the sigma-field $\Cal P^{(-\infty,-m]\cup[m,\infty)}$ 
equals the full (product) sigma-field. Indeed, suppose we see all entries of a point $x$ except on the 
interval $(-m,m)$. In a typical point, this interval is contained 
between a pair of successive markers $i$ for some level $k$ of the inductive construction. Then, by examining 
this point's entries far enough to the left and right we will see completely all but one blocks from the 
family $\Cal B_i$ which constitute the block $C\in\Cal C_i$ covering the considered interval. Because 
every block from $\Cal B_i$ is used in $C$ exactly once, by elimination, we will be able to determine
the missing block and hence all symbols in $(-m,m)$.
\endremark

\medskip
The next construction shows that there exists a process isomorphic to a Bernoulli process with an almost 
sure limit law $\tilde F\equiv 0$ for the normalized hitting times (strong attracting), achieved along a 
subsequence of upper density 1. In particular, this answers in the negative a question of Zaqueu Coelho 
([C]), whether all processes isomorphic to Bernoulli processes have necessarily the exponential limit 
law for the hitting (and return) times. The idea of the construction was suggested to us by D. Rudolph 
([R2]), who attributes the method to Arthur Rothstein. 

\example{Example 2} We will build a decreasing sequence of subshifts of finite type (SFT's). In each we will 
regard the measure of maximal entropy. Begin with the full shift $X_0$ on a finite alphabet. Select $r$
words $W_1, W_2, \dots, W_r$ of some length $l$ and create $r$ SFT's: $X_0^{(1)}, X_0^{(2)},\dots, X_0^{(r)}$, 
forbidding one of these words in each of them, respectively. Choose a length $n$ so large that in the majority 
of blocks of this length in $X_0^{(i)}$ all words $W_j$ except $W_i$ will appear at least once. Now choose
another length $m$, such that in the majority of blocks of this length every block $C$ of length between $n$ 
and $n^2$ will appear many times. Now define $X_1$ as the subshift whose each point is a concatenation of the 
form $\dots B_1B_2\dots B_rB_1\dots$, where $B_i$ is a block appearing in $X_0^{(i)}$ of length either $m$ or 
$m+1$. Obviously, a typical block $C$ of any length between $n$ and $n^2$ appearing in $X_1$ comes from some 
$X_0^{(i)}$, hence contains all $W_j$'s except $W_i$, therefore in a typical $x\in X_1$, $C$ will appear many 
times within each component $B_i$ representing $X_0^{(i)}$, and then it will be absent for a long time, until 
the next representative of $X_0^{(i)}$. So, every such block will reveal strong attracting. It is not hard 
to see that $X_1$ is a mixing SFT and its d-bar distance from the full shift is small whenever the length 
$l$ of the (few) forbidden words $W_i$ is large. We can now repeat the construction starting with $X_1$, and 
radically increasing all parameters. We can arrange that the d-bar distances are summable, so the 
limit system $X$ (the intersection of the $X_k$'s), more precisely its measure of maximal entropy, is also a 
d@-bar limit. Each mixing SFT is isomorphic to a Bernoulli process and this property passes via d-bar limits 
(see [O], [Sh]), hence $X$ is also isomorphic to a Bernoulli process. This system has the almost sure limit law 
$\tilde F\equiv 0$ for hitting times (or $F\equiv 1$ for the return times) achieved along a \sq\ containing infinitely
many intervals of the form $[n,n^2]$. Such \sq\ has upper density 1. 
\endexample

\remark{Remark 3}
It is also possible to construct a process $X_h$ as above with any preassigned entropy $h$. On the other hand, 
it is well known ([Si]), that every measure-preserving transformation with positive entropy $h$ possesses a 
Bernoulli factor of the same entropy. By the Ornstein Theorem ([O]) this factor is isomorphic to 
$X_h$. The generator of $X_h$ appears as a partition of the space on which the initial measure-preserving 
transformation is defined. This proves the universality of ``the law of series'': {\bf in every 
measure-preserving transformation there exists a partition generating the full entropy, which has the
``strong repelling properties''} (i.e., almost sure limit law $\tilde F\equiv 0$ along a \sq\ of lengths of 
upper density 1). 
\endremark

\medskip
Various zero entropy processes with persistent repelling or attracting are implicit in the existing literature.
Extreme repelling (with intensity converging to $e^{-1}$ as the length of blocks grows) occurs for example 
in odometers, or, more generally, in rank one systems ([C-K]). For completeness, we sketch two zero entropy
processes with features of positive entropy: repelling, and the unbiased behavior.

\example{Example 3} Take the product of the independent Bernoulli process on two symbols with an odometer 
(modeled by an adequate process, for example a regular Toeplitz subshift; see [D] for details on Toeplitz flows). 
Call this product process $X_0$. The odometer factor provides, for each $k\in\na$, markers dividing each element 
into so-called {\it $k$-blocks} of equal lengths $p_k$. Each $p_k$ is a multiple of $p_{k-1}$ and each $k$-block 
is a concatenation of $(k\!-\!1)$-blocks. Now we create a new process $X_1$ by ``stuttering'': if 
$x_0\in X_0$ is a concatenation $\dots ABCD\dots$ of $1$-blocks, we create $x_1\in X_1$ as 
$\dots AABBCCDD\dots$, with the number of repetitions 
$q_1=2$. In $X_1$ the lengths of the $k$-blocks for $k>1$ have doubled. Repeating the stuttering for 
$2$-blocks of $X_1$ with a number of repetitions $q_2\ge 2$, we obtain a process $X_2$. And so on. 
Because in each step we reduce the entropy by at least half, the limit process has entropy 
zero. If the $q_k$'s grow sufficiently fast, we obtain, like in the previous example, 
a system with strong attracting for a set of lengths of upper density 1. 
Consider a modification of this example where $q_k = 2$ for each $k$ and each pair $AA$ (also $BB$, 
etc.) is substituted by $A\overline A$ ($B\overline B$, etc.), where $\overline A$ is the ``mirror'' of $A$, 
i.e., with the symbols 0 are replaced by 1 and {\it vice versa}. It is not very hard to compute, that such process 
(although has entropy zero), has the same limit law properties as the independent process: almost sure 
convergence along the full \sq\ $(n)$ to the unbiased (exponential) limit law. 
\endexample

\remark{Remark 4} It is not hard to construct zero entropy processes with persistent mixed behavior. 
For example, applying the ``stuttering technique'' to an odometer one obtains a process in which a typical 
block $B$ occurs in periodically repeated pairs: $BB........BB........BB........$, i.e., with the function
$G_B\approx \min \{1,\frac t2\}$ (which reveals attracting with intensity 
$\frac{\log 2-1}2$ at $t_1=\log2$ and repelling with intensity $e^{-2}$ at $t_2=2$). We skip the details.
\endremark

\heading Questions \endheading
\remark{Question 1} Is there a speed of the convergence to zero of the joint measure of the ``bad'' blocks
in Theorem 1? More precisely, does there exist a positive function $s(n,\epsilon,\#\Cal P)$ converging to 
zero as $n$ grows, such that if for some $\epsilon$ and infinitely many $n$'s, the joint measure of the 
$n$-blocks which repel with intensity $\epsilon$ exceeds $s(n,\epsilon,\#\Cal P)$, then the process has 
necessarily entropy zero? (By the Example 1, $\frac1n$ is not enough.)
\endremark
\remark{Question 2}
Can one strengthen the Theorem 2 as follows: 
$$
\limsup_{n\to\infty} \,\tilde F_{x,n}\le 1-e^{-t} \ \ \mu\text{-almost everywhere?}
$$ 
\endremark
\remark{Question 3}
In Lemma 3, can one obtain $\Cal P^r$ conditionally $\beta$-independent of jointly the past and 
{\it all} return times $R^{(k)}_B$ ($k\ge 1$) (for sufficiently large $n$, with $\mu$-tolerance 
$\beta$ for $B\in\Cal P^{-n}$)? In other words, can the $\beta$-independent process \procbr\ be
obtained $\beta$-independent of the factor-process generated by the partition into $B$ 
and its complement?
\endremark
\remark{Question 4} (suggested by J-P. Thouvenot) Find a purely combinatorial proof of Theorem 1,
by counting the quantity of very long strings (of length $m$) inside which a positive fraction (in measure) 
of all $n$-blocks repel with a fixed intensity. For sufficiently large $n$ this quantity should 
be eventually (as $m\to\infty$) smaller than $h^m$ for any preassigned positive $h$.
\endremark
\remark{Question 5} As we have mentioned, we only know about conditions which ensure that the limit law
for the return time is exponential. It would be interesting to find a (large) class of positive entropy processes 
for which the distributions of return times are essentially deviated from exponential for bounded away from zero 
in measure collections of arbitrarily long blocks, i.e., a class of processes with persistent attracting. 
Can one prove that persistent attracting is, in some reasonable sense, a ``typical'' property in positive 
entropy, or that for a fixed measure-preserving transformation with positive entropy, a ``typical'' generator 
(partition) leads to persistent attracting?
\endremark

\enddocument